\documentclass{gtpart}
\usepackage{graphicx}
\usepackage[mathscr]{eucal} 
\usepackage{amsmath}
\usepackage{caption}
\usepackage{subcaption}
\usepackage{hyperref}
\usepackage{xcolor}
\usepackage{float}
\usepackage{graphicx}
\usepackage{longtable}
\usepackage{multicol}
\usepackage{multirow}

\usepackage{algorithm}
\usepackage{algpseudocode}
\usepackage{algorithmicx}

\usepackage{fontawesome}
\newcommand{\dd}{\mathrm{d}}
\newcommand{\T}{\mathsf{T}}
\newcommand{\RR}[1]{\mathbf{R}^{#1}}
\newcommand{\Cinf}[1]{\mathbf{C}^{\infty}(#1)}
\newcommand{\ec}[1]{\mbox{$#1$}}

\newtheorem{Theorem}{Theorem}
\newtheorem{theorem}[Theorem]{Theorem}

\newtheorem{Proposition}[Theorem]{Proposition}

\newtheorem{example}[Theorem]{Example}
\newtheorem{Remark}{Remark}

\newcommand{\Input}[1]{\mbox{\textbf{Input:}\,\,}\parbox[t]{0.91\linewidth}{#1}}
\newcommand{\Output}[1]{\mbox{\textbf{Output:}\,\,}\parbox[t]{0.88\linewidth}{#1}}
\newcommand{\rreturn}[1]{\textbf{return}\,\,{#1}}

\numberwithin{equation}{section}

\begin{document}
\title{Areas on the space of smooth probability \\ density functions on \ec{\mathbf{S}^2}}
%

\author{J. C. Ru\'iz-Pantale\'on}

\address{Universidad de Sonora}

\email{jose.ruiz\,@unison.mx}

\author{P. Su\'arez-Serrato} 

\address{Instituto de Matem\'aticas, Universidad Nacional Aut\'onoma de M\'exico (UNAM), Mexico City, Mexico}

\email{pablo\,@im.unam.mx}



%
%

%
\begin{abstract}
%
\noindent We present symbolic and numerical methods for computing Poisson brackets on the spaces of measures with positive densities of the plane, the 2--torus, and the 2--sphere. We apply our methods to compute symplectic areas of finite regions for the case of the 2--sphere, including an explicit example for Gaussian measures with positive densities.

\end{abstract}

\maketitle
    \section{Introduction}
    
Since Otto \cite{Otto-2001} introduced a Riemannian structure on the \ec{L^2}--Wasserstein space, plenty of geometric questions have been formulated. A rigorous description of this Riemannian structure was explained by Lott \cite{Lott-2008}, where he also showed that if a manifold $M$ is symplectic or admits a Poisson structure, then its smooth Wasserstein space \ec{P^{\infty}(M)} admits a Poisson bracket, with an associated symplectic foliation. In a related direction, Khesin and Lee study applications of Poisson brackets on density manifolds to geostrophic equations \cite{Khesin-2008}. These works were preceded by ideas of Weinstein \cite{Wein-83}, and Marsden--Weinstein \cite{Marsd-Wein-82}. More details on the symplectic and Poisson structures on Wasserstein spaces may be found in the monograph by Gangbo--Kim--Tomasso \cite{Gangbo-2011} and in the book by Khesin--Wendt \cite[Appendix 5]{Khesin-book}. One motivation comes from applications to geometric hydrodynamics, known already to Arnold \cite{Arnold}. Furthermore, Khesin--Misiolek--Modin discovered profound connections with K\"ahler geometry \cite{Khesin-2018}.

The development of computational aspects of geometric structures on Wasserstein geometry is a very active area. For example, applications to image analysis and signal processing have been extensively developed \cite{CompInfoGeom-book}. An example also includes approximations of Wasserstein geodesics by cubic splines \cite{Benamou-2019}. Furthermore, an algorithm that produces samples with respect to a Hamiltonian flow and a volume form coming from a symplectic form is available \cite{Barp}. Recently, in collaboration with M. A. Evangelista--Alvarado, we have explored the symbolic \cite{EPS1} and numerical \cite{EPS2} computational aspects of Poisson geometry.

In this note, we carry out symbolic and numerical computations of the symplectic form (and associated area) on the space of smooth probability density functions of the 2--sphere. We include the computation of the analogous Poisson bracket for the 2--torus as well. We contribute Algorithms \ref{alg:PWM} and \ref{alg:PWS2}, as well as their symbolic and numerical implementations\footnote{Found in the repository \url{https://github.com/appliedgeometry/density-areas}.}. To the best of our knowledge this is the first time these algorithms have been implemented, and also the first computations of areas of symplectic leaves on inifinite--dimensinoal Poisson manifolds.
In the case of the 2--sphere, we calculate areas in the space of Gaussian measures with positive densities. We exploit the fact that rotations on the 2--sphere are spanned by Hamiltonian flows and therefore preserve measures, in particular, symplectic areas. Our main result in this direction is Theorem \ref{thm:param-area}. As an example, an explicit computation of the area of a patch of a symplectic leaf in \ec{P^{\infty}(S^2)} is found in Equation \ref{eq:area} (see Example \ref{sec:areaGauss-S2}). The performance of numerical computation is generally faster than that of symbolic computation. However, with the latter, we can sometimes obtain precise results, finding exact analytical formul{\ae}, shown in several examples computed in Tables \ref{table:PWR2}, \ref{table:PWT2} and \ref{table:PWS2}.

In section \ref{sec:PW} we briefly review the Poisson geometry of the smooth Wasserstein space and we explain both our Algorithms \ref{alg:PWM} and \ref{alg:PWS2}. Section \ref{sec:areaGauss-S2} includes our Theorem \ref{thm:param-area} which gives the theoretical explanation for calculating areas of symplectic leaves in \ec{P^{\infty}(S^2)}. 




    \section{Poisson Structures on \ec{P^{\infty}(M)}} \label{sec:PW}

Let $M$ be a smooth, connected closed manifold. We denote by \ec{\dd\mathrm{vol}_{M}^{g}} the Riemannian density induced by given a smooth Riemannian metric $g$ on $M$.

The space of absolutely continuous measures with a smooth positive density function on $M$ is defined as:
    \begin{equation*}
        P^{\infty}(M) := \left\{ \rho\,\dd\mathrm{vol}_{M}^{g} \,\Big|\, \rho \in \Cinf{M}, \rho > 0, \int_{M}\rho\,\dd\mathrm{vol}_{M}^{g} = 1 \right\}
    \end{equation*}
It is well--known that \ec{P^{\infty}(M)} is a dense subset of the space of Borel probability measures on $M$, equipped with the Wasserstein metric \ec{W_{2}} \cite{Lott-2008}. Moreover, it has the structure of an infinite--dimensional smooth manifold \cite{Michor-1997}.

\begin{Proposition}
\cite{Marsd-Wein-82,Wein-83,Lott-2008} If $M$ is a Poisson manifold, then \ec{P^{\infty}(M)} admits a Poisson bracket.
\end{Proposition}

Let $\pi$ be a Poisson bivector field  on $M$, and \ec{\{\,,\,\}_{\pi}} the induced Poisson bracket. For \ec{\varphi \in \Cinf{M}} define a linear functional \ec{\mathscr{F}_{\varphi} \in \Cinf{P^{\infty}(M)}} by
    \begin{equation}\label{EcLF}
        \mathscr{F}_{\varphi}(\mu) := \int_{M} \varphi\,\dd{\mu}, \quad \mu \in P^{\infty}(M).
    \end{equation}
Let \ec{f,h \in \Cinf{M}}, then the following formula defines a Poisson bracket on \ec{P^{\infty}(M)}:
    \begin{equation}\label{ec:PoissBracketP00}
        \{\mathscr{F}_{f}, \mathscr{F}_{h}\}(\mu) := \int_{M} \{f,h\}_{\pi}\,\dd{\mu}, \quad \mu \in P^{\infty}(M)
    \end{equation}
By definition, the map \ec{f \to \mathscr{F}_{f}} is Lie algebra morphism, because \ec{\{\mathscr{F}_{f}, \mathscr{F}_{h}\} = \mathscr{F}_{\{f,h\}_{\pi}}}.

    \subsection{Poisson Structures on \ec{P^{\infty}(\RR{2})} and \ec{P^{\infty}(\mathbf{T}^{2})}} \label{subsec:R2}

We will now work on bounded domains on the Euclidean plane and on the 2--torus.

Let \ec{Q = (0,1) \times (0,1)} be the open unit square in \ec{\RR{2}}. Consider the Poisson manifold \ec{(Q, \pi)}, with $\pi$ a Poisson bivector field on $Q$. Without loss of generality assume that our bounded domain of interest lies inside $Q$. Let \ec{(x,y)} be Cartesian coordinates on $Q$. Then $\pi$ has the following representation,
    \begin{equation}\label{EcPiR2}
        \pi = \tau\,\frac{\partial}{\partial x} \wedge \frac{\partial}{\partial y}
    \end{equation}
for some conformal factor \ec{\tau \in \Cinf{Q}}.

Endow $Q$ with the Euclidean metric \ec{g_{0} = \dd{x^{2}} + \dd{y^{2}}}. Then the Riemannian density induced by \ec{g_{0}} is equal to  \ec{|\dd{x} \wedge \dd{y}|}. Therefore,
    \begin{equation*}
        P^{\infty}(Q) = \left\{ \rho\,|\dd{x} \wedge \dd{y}| \ \Big|\ \rho \in \Cinf{Q}, \rho > 0, \int_{0}^{1}\int_{0}^{1}\rho\,\dd{x}\dd{y} = 1 \right\}.
    \end{equation*}
By equation \eqref{ec:PoissBracketP00}, the Poisson bracket induced by $\pi$ on \ec{P^{\infty}(Q)}, evaluated on the linear functionals \ec{\mathscr{F}_{f},\mathscr{F}_{h} \in \Cinf{P^{\infty}(Q)}} at \ec{\mu_{\rho} = \rho\,|\dd{x} \wedge \dd{y}| \in P^{\infty}(Q)}, is given by
    \begin{equation}\label{ec:PoissBracketP00Q}
        \{\mathscr{F}_{f}, \mathscr{F}_{h}\}(\mu_{\rho}) = \int_{0}^{1}\int_{0}^{1} \left( \frac{\partial f}{\partial x} \frac{\partial h}{\partial y} - \frac{\partial f}{\partial y} \frac{\partial h}{\partial x} \right) \tau\rho\,\mathrm{d}x\mathrm{d}y.
    \end{equation}

Now consider a Poisson bivector field $\pi$ on the 2--torus \ec{\mathbf{T}^{2}}. Let \ec{(\theta_{1}, \theta_{2})} be natural coordinates on \ec{\mathbf{T}^{2} \simeq \mathbf{S}^{1} \times \mathbf{S}^{1}}, with \ec{\theta_{1}, \theta_{2} \in \mathbf{R}/2\pi\mathbf{Z}}. Then $\pi$ has a representation analogous to \eqref{EcPiR2} for some conformal factor \ec{\tau \in \Cinf{\mathbf{T}^{2}}}. Hence, endowing \ec{\mathbf{T}^{2}} with the product metric:
    \begin{equation*}
        P^{\infty}(\mathbf{T}^{2}) = \left\{ \rho\,|\dd{\theta_{1}} \wedge \dd{\theta_{2}}| \ \Big|\ \rho \in \Cinf{\mathbf{T}^{2}}, \rho > 0, \frac{1}{4\pi^{2}}\int_{0}^{2\pi}\int_{0}^{2\pi}\rho\,\dd{\theta_{1}}\dd{\theta_{2}} = 1 \right\}
    \end{equation*}
Therefore the Poisson bracket induced by $\pi$ on \ec{P^{\infty}(\mathbf{T}^{2})} is given by:
    \begin{equation}\label{ec:PoissBracketP00T2}
        \{\mathscr{F}_{f}, \mathscr{F}_{h}\}(\mu_{\rho}) = \frac{1}{4\pi^{2}}\int_{0}^{2\pi}\int_{0}^{2\pi} \left( \frac{\partial f}{\partial \theta_{1}} \frac{\partial h}{\partial \theta_{2}} - \frac{\partial f}{\partial \theta_{2}} \frac{\partial h}{\partial \theta_{1}} \right) \tau\rho\,\mathrm{d}\theta_1\mathrm{d}\theta_2
    \end{equation}

Our next algorithm computes the Poisson bracket in \eqref{ec:PoissBracketP00Q} and \eqref{ec:PoissBracketP00T2}.
\begin{algorithm}[H]\label{algo1}
    \captionsetup{justification=centering}
    \caption{\ \textsf{2DEuclidean\_PoissonWasserstein}} \label{alg:PWM}
        \rule{\textwidth}{0.4pt}
    \Input{$\tau$ conformal factor, $\rho$ positive density function, \ec{f,h} scalar functions on \ec{Q,\mathbf{T}^{2}}}
    \Output{the value of the Poisson bracket on either \ec{P^{\infty}(Q)} or \ec{P^{\infty}(\mathbf{T}^{2})}, induced by a Poisson bivector field \ec{\pi_{\tau}} on $Q$ or \ec{\mathbf{T}^{2}}, of \ec{\mathscr{F}_{f}} and \ec{\mathscr{F}_{h}} in \ec{\Cinf{P^{\infty}(Q)}, \Cinf{P^{\infty}(\mathbf{T}^{2})}} at \ec{\mu_{\rho}} in \ec{P^{\infty}(Q), P^{\infty}(\mathbf{T}^{2})}, in that order}
        \rule{\textwidth}{0.4pt}
    \begin{algorithmic}[1] 
        \Procedure{}{}
            \State \textsc{Transform} \ec{\tau, \rho, f} and $h$ to symbolic variables
            \If {\ec{\mathbf{T}^{2}}}
                \State \textsc{Transform} \ec{\theta_{1}} and \ec{\theta_{2}} to symbolic variables
                \State \emph{integrand} $\leftarrow$ \ec{\big( \frac{\partial f}{\partial \theta_{1}} \frac{\partial h}{\partial \theta_{2}} - \frac{\partial f}{\partial \theta_{2}} \frac{\partial h}{\partial \theta_{1}} \big)\tau\rho}
                \If {numerical=True}
                    \State \rreturn{a numerical approximation of the double integral \ec{\frac{1}{4\pi^{2}}\int_{0}^{2\pi}\int_{0}^{2\pi} \text{\emph{integrand}}\,\mathrm{d}\theta_1\mathrm{d}\theta_2}}
                \EndIf
                \State \rreturn{the double integral \ec{\frac{1}{4\pi^{2}}\int_{0}^{2\pi}\int_{0}^{2\pi}\text{\emph{integrand}}\,\mathrm{d}\theta_1\mathrm{d}\theta_2}}
            \EndIf
            \State \textsc{Transform} \ec{x} and \ec{y} to symbolic variables
            \State \emph{integrand} $\leftarrow$ \ec{\big( \frac{\partial f}{\partial x} \frac{\partial h}{\partial y} - \frac{\partial f}{\partial y} \frac{\partial h}{\partial x} \big)\tau\rho}
            \If {numerical=True}
                \State \rreturn{a numerical approximation of the double integral \ec{\int_{0}^{1}\int_{0}^{1} \text{\emph{integrand}}\,\mathrm{d}x\mathrm{d}y}}
            \EndIf
            \State \rreturn{the double integral \ec{\int_{0}^{1}\int_{0}^{1} \text{\emph{integrand}}\,\mathrm{d}x\mathrm{d}y}}
        \EndProcedure
    \end{algorithmic}
\end{algorithm}

\begin{example}
For \ec{(Q, \pi = \frac{\partial}{\partial x} \wedge \frac{\partial}{\partial y})} consider \ec{\rho=\tfrac{3}{2}(x^{2} + y^{2})}, a radial density function in Cartesian coordinates, and \ec{\rho=\frac{2}{c\pi}\mathrm{exp}[-\frac{1}{2}(x^{2} + y^{2})]}, a bivariate normal density with \ec{c = \mathrm{Erf}(\sqrt{2}/{2})^{2}} and $\mathrm{Erf}$ equal to Gauss' error function. The factors \ec{3/2} and \ec{1/c} are normalisation constants that make these densities have total integral equal to 1 over $Q$. This is the input data in Table \ref{table:PWR2}, where we display symbolic and numerical values for \ec{\{\mathscr{F}_{f}, \mathscr{F}_{h}\}(\mu_{\rho})}, with respect to some simple functions \ec{f,g}.
    \begin{table}[H]
        \centering
        \caption{Symbolic and numeric computational examples of the Poisson bracket \eqref{ec:PoissBracketP00Q} with \ec{\tau = 1}; for densities $\rho$ and scalar functions $f$ and $h$, all on \ec{\mathbf{R}^{2}}. } \label{table:PWR2}
    \begin{tabular}{|c|c|c|c|c|}
            \cline{4-5}
        \multicolumn{1}{c}{} & \multicolumn{1}{c}{} & \multicolumn{1}{c}{} & \multicolumn{2}{|c|}{\ec{\{\mathscr{F}_{f}, \mathscr{F}_{h}\}(\mu_{\rho})}} \\
            \hline
        \multicolumn{1}{|c|}{$\rho$} & \multicolumn{1}{c}{$f$} & \multicolumn{1}{|c|}{$h$} & Symbolic & Numerical (estimate, error) \\
            \hline
            \hline
        1 & $x$ & $y$ & $\phantom{-}1$ & $(1.0, 1.1102230246251565e^{-14})$ \\
             \hline
        $\tfrac{3}{2}(x^{2} + y^{2})$ & $x + y$ & $x^2 - y^2$ & $-\frac{5}{2}$ & $(-2.5000000000004,6.90760897044e^{-14})$ \\
            \hline
        '' & $\sin{x}$ & $\sin{y}$ & $3(\sin{2} - \sin^{2}{1})$ & $(\phantom{-}0.6036720256563,9.72301012762e^{-15})$ \\
            \hline
        $\frac{2}{c\pi}e^{-\frac{1}{2}(x^{2} + y^{2})}$ & $\sin{x}$ & $\cos{x}$ & $\phantom{-}0$ & $(0.0,0)$\\
            \hline
        '' & $x^{2}$ & $y^{2}$ & $\frac{8}{c\pi}\frac{(\sqrt{e} - 1)^{2}}{e}$ & $(\phantom{-}0.8458930796971,1.44765646605e^{-14})$\\
            \hline
    \end{tabular}
    \end{table}
\end{example}

\begin{example}
We can verify that for \ec{\tau = 1} and linear functions \ec{f=ax + by} and \ec{h=cx+dy} on $Q$, with \ec{a,b,c,d \in \mathbf{R}}, the Poisson bracket of \ec{\mathscr{F}_{f}} and \ec{\mathscr{F}_{h}} in \ec{\mathbf{C}^{\infty}(P^{\infty}(Q))} at an arbitrary \ec{\mu_{\rho}} in \ec{P^{\infty}(Q)} is given by \ec{\{\mathscr{F}_{f}, \mathscr{F}_{h}\}(\mu_{\rho}) = ad - bc}.
\end{example}

\begin{example}
For \ec{(\mathbf{T}^{2}, \pi = \frac{\partial}{\partial \theta_{1}} \wedge \frac{\partial}{\partial \theta_{2}})} consider \ec{\rho = \tfrac{1}{2} (\cos\theta_1 + \cos\theta_2 + 2 )}, a density that equals the sum of two cardioid density functions on \ec{\mathbf{S}^{1}} \cite{Jeffreys}, and \ec{\rho = \frac{1}{c}e^{\cos\theta_1 + \cos\theta_2}}, a (particular case of a) density featured in probabilistic modelling of torsional angles in molecules with $c = 1/I_{0}(1)^2 \approx 0.6238603604320$ and \ec{I_{0}} equal to the (0--th) modified Bessel function of the first kind \cite{Singh}. The factors \ec{1/2} and \ec{1/c} are normalisation constants that make these densities have total integral equal to 1 over \ec{\mathbf{T}^{2}}. This is the input data in Table \ref{table:PWT2}, where we display symbolic and numerical values for \ec{\{\mathscr{F}_{f}, \mathscr{F}_{h}\}(\mu_{\rho})}, with respect to some simple functions \ec{f,g}.
    \begin{table}[H]
        \centering
        \caption{Symbolic and numeric computational examples of the Poisson bracket \eqref{ec:PoissBracketP00T2} with \ec{\tau = 1}; for densities $\rho$ and scalar functions $f$ and $h$, all on \ec{\mathbf{T}^{2}}.} \label{table:PWT2}
     \begin{tabular}{|c|c|c|c|c|}
            \cline{4-5}
        \multicolumn{1}{c}{} & \multicolumn{1}{c}{} & \multicolumn{1}{c}{} & \multicolumn{2}{|c|}{\ec{\{\mathscr{F}_{f}, \mathscr{F}_{h}\}(\mu_{\rho})}} \\ \cline{1-5}
        \multicolumn{1}{|c|}{$\rho$} & \multicolumn{1}{c}{$f$} & \multicolumn{1}{|c|}{$h$} & Symbolic & Numerical (estimate, error) \\
            \hline
            \hline
        1 & $\theta_{1}$ & $\theta_{2}$  & $\phantom{-}1$ & $(1.0,1.11022302462e^{-14})$  \\
             \hline
        $\scriptstyle\tfrac{1}{2} ( \cos\theta_1 + \cos\theta_2 + 2 )$ & $\theta_{1}^{2}$ & $\theta_{2}^{2}$  & $4\pi^{2}$ & $(39.47841, 4.382984e^{-13})$ \\
            \hline
        '' & $e^{\theta_{1}}$ & $e^{\theta_{2}}$  & $\frac{3}{8\pi^2}(e^{2\pi} - 1)^2$ & $(10854.5889,1.3894e^{-10})$ \\
            \hline
        $\frac{1}{c}e^{\cos\theta_1 + \cos\theta_2}$ & $\theta_{1}$ & $\theta_{2}$  & $1$ & $(1.0,3.062998047371e^{-12})$ \\
            \hline
        '' & $\theta_{1}^{2}$ & $\theta_{2}^{2}$  & $\frac{c}{\pi^{2}}\int_{0}^{2\pi}\theta_{1}e^{\cos{\theta_{1}}} \dd{\theta_{1}} \int_{0}^{2\pi}\theta_{2}e^{\cos{\theta_{2}}} \dd{\theta_{2}} $ & $(39.47841, 2.222763e^{-10})$ \\
            \hline
    \end{tabular}
    \end{table}
\end{example}

    \subsection{Poisson Structures on \ec{P^{\infty}(\mathbf{S}^{2})}} \label{subsec:S2}

Let \ec{(\theta, \phi) \in (0,\pi) \times (0,2\pi)} be spherical coordinates on \ec{\mathbf{S}^{2}}, such that
\[
(\theta, \phi) \mapsto ( \sin{\theta} \cos{\phi}, \sin{\theta} \sin{\phi}, \cos{\theta} ).
\]
Consider the round metric \ec{g_{0} = \dd{\theta}^{2} + \sin^{2}{\theta}\,\dd{\phi}^{2}} on \ec{\mathbf{S}^{2}} seen as a unit sphere. The metric \ec{g_{0}} induces the Riemannian area form
    \begin{equation}\label{eqn:can-symp-s2}
        \dd\mathrm{vol}_{\mathbf{S}^{2}}^{g_{0}} = \sin{\theta}\,\dd{\theta} \wedge \dd{\phi}.
    \end{equation}
 Denote by \ec{\pi_{0}} the Poisson bivector induced by the symplectic form (\ref{eqn:can-symp-s2}) on \ec{\mathbf{S}^{2}}.
Let \ec{\tau \in \Cinf{\mathbf{S}^{2}}} be a conformal factor, then \ec{(\mathbf{S}^{2}, \pi_{\tau} = \tau \pi_{0})} is a Poisson manifold. Therefore \ec{\pi_{\tau}=\frac{\tau}{\sin{\theta}}\,\frac{\partial}{\partial \theta} \wedge \frac{\partial}{\partial \phi}}, and we may write:
    \begin{equation*}
        P^{\infty}\left( \mathbf{S}^{2} \right) = \left\{ \rho\sin{\theta}\,|\dd{\theta} \wedge \dd{\phi}| \ \Big|\ \rho \in \Cinf{\mathbf{S}^{2}}, \rho > 0, \frac{1}{4\pi}\int_{0}^{2\pi}\int_{0}^{\pi}\rho\sin{\theta}\,\dd{\theta}\dd{\phi} = 1 \right\}
    \end{equation*}

Let \ec{f,h \in \Cinf{\mathbf{S}^{2}}}. By equation \eqref{ec:PoissBracketP00}, the Poisson bracket induced by \ec{\pi_{\tau}} on \ec{P^{\infty}(\mathbf{S}^{2})} and evaluated on \ec{\mathscr{F}_{f},\mathscr{F}_{h} \in \Cinf{P^{\infty}(\mathbf{S}^{2})}}
at \ec{\mu_{\rho} = \rho\,\dd\mathrm{vol}_{\mathbf{S}^{2}}^{g_{0}} \in P^{\infty}(\mathbf{S}^{2})} is given by:
    \begin{equation}\label{ec:PoissBracketP00S2}
        \{\mathscr{F}_{f}, \mathscr{F}_{h}\}(\mu_{\rho}) = \frac{1}{4\pi}\int_{0}^{2\pi}\int_{0}^{\pi} \left( \frac{\partial f}{\partial \theta} \frac{\partial h}{\partial \phi} - \frac{\partial f}{\partial \phi} \frac{\partial h}{\partial \theta} \right) \tau\rho\,\mathrm{d}\theta\mathrm{d}\phi
    \end{equation}

Our following algorithm computes the Poisson bracket expressed in \eqref{ec:PoissBracketP00S2}.
\begin{algorithm}[H]\label{algo2}
    \captionsetup{justification=centering}
    \caption{\ \textsf{2DSpherical\_PoissonWasserstein}} \label{alg:PWS2}
        \rule{\textwidth}{0.4pt}
    \Input{a conformal factor $\tau$, a positive density function $\rho$, and two scalar functions $f$ and $g$ on \ec{\mathbf{S}^{2}}}
    \Output{the value given in \eqref{ec:PoissBracketP00S2} of the Poisson bracket induced by \ec{\pi_{\tau}} on \ec{P^{\infty}(\mathbf{S}^{2})} of \ec{\mathscr{F}_{f},\mathscr{F}_{h} \in \Cinf{P^{\infty}(\mathbf{S}^{2})}} at \ec{\mu_{\rho} = \rho\,\dd\mathrm{vol}_{\mathbf{S}^{2}}^{g_{0}} \in P^{\infty}(\mathbf{S}^{2})}}
        \rule{\textwidth}{0.4pt}
    \begin{algorithmic}[1] 
        \Procedure{}{}
            \State \textsc{Transform} \ec{\theta, \phi, \tau, \rho, f} and $h$ to symbolic variables
            \State \emph{integrand} $\leftarrow$ \ec{\big( \frac{\partial f}{\partial \theta} \frac{\partial h}{\partial \phi} - \frac{\partial f}{\partial \phi} \frac{\partial h}{\partial \theta} \big)\tau\rho}
            \If {numerical=True}
                \State \rreturn{a numerical approximation of the double integral \ec{\frac{1}{4\pi}\int_{0}^{2\pi}\int_{0}^{\pi} \text{\emph{integrand}}\,\mathrm{d}\theta\mathrm{d}\phi}}
            \EndIf
            \State \rreturn{the double integral \ec{\frac{1}{4\pi}\int_{0}^{2\pi}\int_{0}^{\pi} \text{\emph{integrand}}\,\mathrm{d}\theta\mathrm{d}\phi}}
        \EndProcedure
    \end{algorithmic}
\end{algorithm}

\begin{example}
For \ec{(\mathbf{S}^{2}, \pi_{0} = \frac{1}{\sin{\theta}}\frac{\partial}{\partial \theta} \wedge \frac{\partial}{\partial \phi})} consider a radial density in polar coordinates \ec{\frac{2}{\pi}\frac{\sin{\theta}}{1-\cos{\theta}}}, and an exponential density \ec{\frac{1}{c\pi}e^{\frac{\sin{\theta}}{\cos{\theta}-1}}}, both on \ec{\mathbf{S}^{2}}, with \ec{c \approx 0.19781355915936985}. The previous two densities are obtained using the stereographic projection starting from densities that appear in the study of isoperimetric regions in the plane \cite{Morgan}. The factors \ec{2/\pi} and \ec{1/c\pi} are normalisation constants that make these densities have total integral equal to 1 over \ec{\mathbf{S}^{2}}. Table \ref{table:PWS2} displays symbolic and numerical values for \ec{\{\mathscr{F}_{f}, \mathscr{F}_{h}\}(\mu_{\rho})}, with respect to some simple functions \ec{f,h}.
    \begin{table}[H]
        \small
        \centering
        \caption{Symbolic and numeric computational examples of the Poisson bracket \eqref{ec:PoissBracketP00S2} with \ec{\tau = 1}; for densities $\rho$ and scalar functions $f$ and $h$, all on \ec{\mathbf{S}^{2}}. Notice that the symbolic value in the third row did not converge.} \label{table:PWS2}
     \begin{tabular}{|c|c|c|c|c|}
            \cline{4-5}
        \multicolumn{1}{c}{} & \multicolumn{1}{c}{} & \multicolumn{1}{c}{} & \multicolumn{2}{|c|}{\ec{\{\mathscr{F}_{f}, \mathscr{F}_{h}\}(\mu_{\rho})}} \\ \cline{1-5}
        \multicolumn{1}{|c|}{$\rho$} & \multicolumn{1}{c}{$f$} & \multicolumn{1}{|c|}{$h$} & Symbolic & Numerical (estimate, error) \\
            \hline
            \hline
        1 & $\theta$ & $\phi$  & $\frac{\pi}{2}$ & $(1.5707963, 1.7439342e^{-14})$  \\
             \hline
        $\frac{2}{\pi}\frac{\sin{\theta}}{1-\cos{\theta}}$ & $\sqrt{\frac{5}{16\pi}}\,(3\cos^{2}\theta - 1)$ &  $\sqrt{\frac{15}{16\pi}}\,\cos(2x)\sin^{2}\theta$  & $0$ & $(-3.2102291e^{-18}, 1.7467849e^{-11})$ \\
            \hline
        $\frac{1}{c\pi}e^{\frac{\sin{\theta}}{\cos{\theta}-1}}$ & $\sqrt{\frac{3}{4\pi}}\sin{\theta}\sin{\phi}$ & $\sqrt{\frac{3}{4\pi}}\sin{\theta}\cos{\phi}$ & NA & $(0.0534111, 1.6277804e^{-10})$ \\
            \hline
    \end{tabular}
    \end{table}
\end{example}

    \section{Areas of smooth measures on the 2-sphere} \label{sec:areaGauss-S2}

Consider \ec{(\mathbf{S}^{2}, \pi_{0} = \frac{1}{\sin{\theta}}\frac{\partial}{\partial \theta} \wedge \frac{\partial}{\partial \phi})} and a measure \ec{\mu_{\rho}} in \ec{P^{\infty}(\mathbf{S}^{2})}. Let the maps
    \begin{equation*}
        T_{s},R_{t}: \mathbf{S}^{2} \to \mathbf{S}^{2}
    \end{equation*}
be spanned by flows of Hamiltonian vector fields on $\mathbf{S}^{2}$ with \ec{s \in I_{1} \subset (0,\pi)}, \ec{t \in I_{2} \subset (0,2\pi)}. The orbits of the natural action on \ec{P^{\infty}(\mathbf{S}^{2})} of the group of Hamiltonian diffeomorphism on \ec{\mathbf{S}^{2}} are the symplectic leaves of \eqref{ec:PoissBracketP00S2}. Therefore, the actions of \ec{T_{s}} and \ec{R_{t}} on \ec{\mu_{\rho}} trace an \ec{(s,t)}--parametrized, 2--dimensional, region $\Omega$ on a symplectic leaf of the Poisson structure on \ec{P^{\infty}(\mathbf{S}^{2})} induced by \ec{\pi_{0}}.

\begin{theorem}\label{thm:param-area}
The symplectic area  \ec{A_{\mu_{\rho}}(\Omega)} of $\Omega$ in a symplectic leaf that contains \ec{\mu_{\rho}} can be calculated as follows:
    \begin{align}\label{EcAreaPooS2}
        A_{\mu_{\rho}}(\Omega) &= \int_{I_{1} \times I_{2}}\{\mathscr{F}_{\theta}, \mathscr{F}_{\phi}\}\big( (T_{s} \circ R_{t})^{\ast}\mu_{\rho} \big) \dd{s}\dd{t}  \\
        &= \frac{1}{4\pi}\int_{I_{1}}\int_{I_{2}}\int_{0}^{2\pi}\int_{0}^{\pi} (T_{s} \circ R_{t})^{\ast}\rho \ \dd{\theta}\dd{\phi}\dd{s}\dd{t} \nonumber
    \end{align}
\end{theorem}
\begin{proof}
We use the description of symplectic leaves and their symplectic forms given by J. Lott in \cite[Proposition 6]{Lott-2008}:
as the maps \ec{T_{s}} and \ec{R_{t}} are spanned by Hamiltonian flows, they preserve the symplectic form in \eqref{eqn:can-symp-s2} and are Hamiltonian diffeomorphisms of \ec{\mathbf{S}^{2}}. Hence, the (finite) region \ec{\Omega := (T_{s} \circ R_{t})^{\ast}\mu_{\rho}} is contained in the symplectic leaf \ec{S_{\mu_{\rho}}} through \ec{\mu_{\rho}} with respect to the Poisson bracket in \eqref{ec:PoissBracketP00S2} with \ec{\tau = 1}. Moreover, the symplectic form on \ec{S_{\mu_{\rho}}} is given by
    \begin{equation*}
        \bar{\omega}(\bar{X}_{f}, \bar{X}_{h}) = \{\mathscr{F}_{f},\mathscr{F}_{h}\}(\mu_{\rho}),
    \end{equation*}
Here \ec{\bar{X}_{f}, \bar{X}_{h} \in \T_{\mu_{\rho}}P^{\infty}(\mathbf{S}^{2})} are infinitesimal motions of \ec{\mu_{\rho}} under the flows generated by the Hamiltonian vector fields \ec{X_{f}} and \ec{X_{h}} on $M$. Therefore the integral in \eqref{EcAreaPooS2} computes the symplectic area of $\Omega$ with respect to $\bar{\omega}$, as claimed.
\end{proof}
\begin{example}[\emph{Areas of smooth Gaussian measures on the 2-sphere.}]\label{Ex::7}
On \ec{(\mathbf{S}^{2}, \pi_{0} = \frac{1}{\sin{\theta}}\frac{\partial}{\partial \theta} \wedge \frac{\partial}{\partial \phi})} consider the following Gaussian measure
    \begin{equation*}
        \mu_{\rho} = \frac{1}{c}\exp\left[ -\frac{1}{4}\left( \frac{\sin\theta}{1-\cos\theta} \right)^{4}\sin^{2}(2\phi) \right]\sin{\theta}\,|\dd{\theta} \wedge \dd{\phi}| \in P^{\infty}(\mathbf{S}^{2}),
    \end{equation*}
with \ec{c \approx 0.7082398710278981}. Define the maps
\[
\ec{T_{s}: \theta \to \arccos(\theta + s)}\quad {\rm and}\quad \ec{R_{t}: \phi \to \phi + t},
\]
with \ec{s,t \in [\pi/4,3\pi/4]}. Set \ec{\Omega = T_s([\pi/4,3\pi/4]) \times R_t([\pi/4,3\pi/4])}. Our Theorem \ref{thm:param-area} and its implementation yield
    \begin{equation}\label{eq:area}
        A_{\mu_{\rho}}(\Omega) \approx 3.3009295509955767,
    \end{equation}
with an estimated error of \ec{5.6961722425556025e^{-9}}.
\end{example}

\begin{Remark}
These numerical experiments ran on a cloud server provided by Google Colaboratory equipped with 12.72 GB of main memory in an Intel(R) Xeon(R) CPU {\makeatletter @} 2.30GHz. Our symbolic implementation is also available. However, it does not seem to converge for the data in equation  \ref{eq:area} of Example \ref{Ex::7} above. This example further illustrates the advantage of numerical methods. The execution time of the computation for the area of \eqref{eq:area} was 4 hours, 27 minutes and 51 seconds.
\end{Remark}

\subsubsection*{Acknowledgments}

This research was supported in part by  DGAPA-UNAM PAPIIT grant IN104819.


\end{document}